\documentclass[10pt, twoside]{article}
\usepackage{hyperref,amssymb,amsmath,graphicx,verbatim,amsfonts,amsthm}
\usepackage{amsmath}
\usepackage{amsfonts}
\newtheorem{thm}{Theorem}%[section]

\newtheorem{lem}[thm]{Lemma}
\newtheorem{prop}[thm]{Proposition}

\newtheorem*{thm*}{Theorem}
\newtheorem*{rem*}{Remark}

\theoremstyle{definition}

\theoremstyle{remark}

\newcommand{\Claim}{\textbf{Claim.}\hspace{10pt}}
\newcommand{\Rem}{\textbf{Remark.}\hspace{10pt}}

\newcommand{\abs}[1]{\left\vert#1\right\vert}

\newcommand{\set}[1]{\left\{#1\right\}}
\newcommand{\Set}[2]{ \left\{#1 \ \big| \ #2 \right\} }
\newcommand{\br}[1]{\left[#1\right]}

\newcommand{\sr}[1]{\left(#1\right)}

\newcommand{\Integer}{\mathbb{Z}}

\newcommand{\Z}{\Integer}
\newcommand{\N}{\mathbb{N}}

\newcommand{\R}{\mathbb{R}}

\newcommand{\eps}{\varepsilon}

\newcommand{\eqdef}{\stackrel{\mathrm{def}}{=}}

\DeclareMathOperator*{\E}{\mathbb{E}}

\renewcommand{\Pr}{}
\let\Pr\relax
\DeclareMathOperator*{\Pr}{\mathbb{P}}

\newcommand{\1}[1]{\mathbf{1}_{\set{ #1 } }}

\def\squareforqed{\hbox{\rlap{$\sqcap$}$\sqcup$}}
\def\qed{\ifmmode\squareforqed\else{\unskip\nobreak\hfil
\penalty50\hskip1em\null\nobreak\hfil\squareforqed
\parfillskip=0pt\finalhyphendemerits=0\endgraf}\fi}

\newcommand{\D}{\mathcal{D}}
\newcommand{\x}{\xi}
\renewcommand{\S}{\mathcal{S}}
\newcommand{\T}{\mathcal{T}}
\newcommand{\U}{\mathcal{U}}

\begin{document}
% \author{Ariel Yadin}

% =====
%\input{paperECP-07-89-head.sty}
% =====

\title{Maximal Arithmetic Progressions in Random Subsets}

\author{Itai Benjamini\footnote{Faculty of Mathematics and Computer Science,
The Weizmann Institute of Science,  POB 26, Rehovot 76100, ISRAEL;
\texttt{itai.benjamini@weizmann.ac.il} }  \and Ariel Yadin
\footnote{Faculty of Mathematics and Computer Science, The
Weizmann Institute of Science,  POB 26, Rehovot 76100, ISRAEL;
\texttt{ariel.yadin@weizmann.ac.il}}
 \and Ofer Zeitouni\footnote{
Department of Mathematics, University of Minnesota, MN 55455, USA;
\texttt{zeitouni@math.umn.edu} } \thanks{Partially supported by
NSF grant DMS-050377 } }

\date{ }

\maketitle

\begin{abstract}
    Let $U^{(N)}$ denote the   maximal length
    of arithmetic progressions in a random
uniform subset of $\set{0,1}^N$.   By an application
of the Chen-Stein method, we show that
%$U^{(N)}/\log N$ converges to $2/\log 2$, whereas
$U^{(N)} -2\log N/\log 2$ converges in law to an extreme type
(asymmetric) distribution. The same result holds for the maximal
length $W^{(N)}$ of arithmetic progressions $\pmod N$. When
considered in the natural way on a common probability space, we
observe that $U^{(N)}/\log N$ converges almost surely to $2/\log
2$, while $W^{(N)}/\log N$ does not converge almost surely (and in
particular, $\limsup W^{(N)}/\log N\geq 3/\log 2$).
\end{abstract}

\section{Introduction and Statement of Results }

In this note we study the length of maximal arithmetic
progressions in a random uniform subset of $\{0,1\}^N$.
%and in
%a random uniform subset of the one dimensional
%torus ($\Z/N\Z$).
That is, let $\x_1,\x_2,\ldots,\x_N$ be a random
word in $\set{0,1}^N$, chosen uniformly. Consider the (random) set $\Xi_N$ of
elements $i$
%in
%$\set{0,1}^N$
%, chosen uniformly. Consider the set of
%$\Z/N\Z$
such that $\x_i = 1$. Let $U^{(N)}$ denote the
%, and find the
maximal length arithmetic progression in $\Xi_N$, and let
$W^{(N)}$
denote the
maximal length aperiodic arithmetic progression $\pmod N$
in $\Xi_N$. A  consequence of our main result
(Theorem \ref{thm: limit thm})
%Proposition
%\ref{prop: LLN})
is that the expectation of both $U^{(N)}$ and
$W^{(N)}$ is
%be the length of this progression. It turns out (see Proposition
%\ref{prop: LLN}) that the expectation of $W^{(N)}$ is
roughly ${2\log N}/{\log 2} $, twice the expectation of the longest
run in $\Xi_N$, see \cite{er},\cite{er1}. We also show
that the limit law of
the centered version of both $W^{(N)}$ and $U^{(N)}$ is of the same
extreme type as that of the longest run in $\Xi_N$.

We observe two interesting phenomena:

\begin{itemize}

\item Theorem \ref{thm: limit thm} states that the tails of the
distribution of $W^{(N)}$ behave differently for positive and
negative deviations from the mean.  In particular, the probability
that $W^{(N)}$ deviates by $x$ from its mean, behaves roughly like
$1 - \exp \sr{ - 2^{-(x+2)} }$ for positive $x$, and like $\exp
\sr{ - 2^{-(x+2)} }$ for negative $x$.  Thus, on the positive side
of the mean the tail decays exponentially, and on the negative
side of the mean the tail decays doubly-exponential.

\item One may construct the sets $\Xi_N$ on the same probability
    space by considering an infinite sequence of
    i.i.d., Bernoulli random variables $\{\xi_i\}_{i=1}^\infty$.
    Proposition
\ref{prop:  Strong LLN for U} states that with such a
construction, the sequence
%OO
${W^{(N)}}/{ \log N}$ converges in
\emph{probability} to the constant ${2}/{\log 2}$, but a.s.
convergence does not hold.  This contrasts with the behavior of
$U^{(N)}$, where a.s. convergence of
%OO
${U^{(N)}}/{ \log N}$
 to  ${2}/{\log 2}$
holds.
%what happens
%when taking arithmetic progressions of subsets on the interval
%(and not the torus).  On the interval, a.s. convergence does hold.
The seemingly small change of taking arithmetic progressions that
``wrap around'' the torus, changes the behavior of the $\limsup$
of the sequence.

\end{itemize}

The notoriously hard extremal problem, showing that a set of
integers of upper positive density contains unbounded arithmetic
progressions, and its finite quantitative versions, is a well
studied topic reviewed recently in \cite{tao}.

\subsection{The Model}

Let $\x_1, \ldots, \x_N, \ldots,$ be i.i.d. Bernoulli random
variables of mean $\E \br{\x_i} = \frac{1}{2}$.   For a
non-negative integer $N$ and $s,p \in \set{1,2,\ldots,N}$ define
$$ W_{s,p} = W_{s,p}^{(N)} \eqdef
\max \Set{ 1 \leq k \leq N }{ \x_s = 0 \ , \
\prod_{i=1}^k \x_{s+i p {\pmod N} } = 1 } . $$ %
That is, we consider all arithmetic progressions $\pmod N$ in
$\set{1,2,\ldots,N}$ starting at $s$, with difference $p$, and
check for the longest one of the form $0,1,1,\ldots$ (the role of
the $0$ is to avoid considering periodic progressions).  $W_{s,p}$
is the length of such progression. We set
$$  W^{(N)} \eqdef \max_{s, p} W_{s,p} , $$
which is the size of the maximal arithmetic progression $\pmod N$ in
$\set{1,2,\ldots,N}$ of the form $0,1,1,\ldots$.

\begin{rem*}
%\label{rem:restricting p}
If $s,s+p,\ldots,s+kp$ is an arithmetic progression ${\pmod N}$ such that
$$ (\x_s,\x_{s+p {\pmod N}} , \ldots, \x_{s+ k p {\pmod N}} , \x_{s+(k+1)p {\pmod N}} ) = (0,1,\ldots,1,0) , $$
then we can also consider the arithmetic progression ${\pmod N}$: starting at $s':=s+(k+1)p {\pmod N}$,
and continuing by jumps of $p' := N-p$.  So $s' + i p' \equiv s + (k+1-i) p {\pmod N}$, and
\begin{align*}
& (\x_{s'} , \x_{s'+p' {\pmod N} } , \ldots, \x_{s' + k p' {\pmod N}} , \x_{s' + (k+1) p' {\pmod N}} ) \\
= 
& 
(\x_{s+(k+1)p {\pmod N}} , \x_{s+kp {\pmod N}} , \ldots, \x_{s+p {\pmod N}} , \x_s ) = (0,1,\ldots,1,0) .
\end{align*}
Thus, since the maximal arithmetic progression ${\pmod N}$ must be followed by a $0$, this progression can possibly be given by two
different jumps $p$, one which moves through the progression in one direction, and the other
in the other direction.  To avoid this, we may consider only jumps $p$ that are $p \leq N/2$.
This restriction will still capture the maximal arithmetic  progression ${\pmod N}$, without going over  the maximal progression twice.
That is, we have that
$$ W^{(N)} = \max_{\substack{ 1 \leq s \leq N \\ 1 \leq p \leq N/2 } } W_{s,p} . $$
Throughout, we will restrict to $p\leq N/2$ in our analysis, the reader should keep in mind this remark.
\end{rem*}

Similarly,
define
%OO
$$ U_{s,p} = U_{s,p}^{(N)} \eqdef \max \Set{ 1 \leq k \leq \lfloor
\frac{N-s}{p}\rfloor}{ \x_s = 0 \ , \
\prod_{i=1}^k \x_{s+i p} = 1 } , $$ %
and
%OO
$$ U = U^{(N)} \eqdef \max_{s, p\leq N} U_{s,p} , $$
that is we only consider $s,p,k$ such that
$\Set{s+ip}{i=0,1,\ldots,k} \subseteq \set{1,\ldots,N}$.

\subsection{Results}

Throughout, we set $C=2/\log 2$.
Our first main result is the following extreme type limit theorem.
\begin{thm} \label{thm: limit thm}
 Let $\set{x_N}$ be a sequence such
that $C \log N + x_N \in \Z$ for all $N$, and $\inf_N x_N \geq b$,
for some $b \in \R$. Then, with
$\lambda(x) = 2^{- (x+3)}$, we have
\begin{equation}\label{elal1}
    \lim_{N \to \infty}
\exp \sr{\lambda(x_N) } \Pr \br{ W^{(N)} \leq C \log N + x_N } = 1\,.
\end{equation}
Similarly, let
 $\set{y_N}$ be a sequence such
that $C \log N -\log_2(2C\log N)+ y_N \in \Z$
for all $N$, and $\inf_N y_N \geq b$,
for some $b \in \R$. Then,
\begin{equation}\label{elal2}
    \lim_{N \to \infty}
\exp \sr{\lambda(y_N) } \Pr \br{ U^{(N)} \leq C \log N -
\log_2(2C\log N)+ y_N } = 1\,.
\end{equation}
In particular, both $W^{(N)}/\log N$ and $U^{(N)}/\log N$
converge in probability
to $C$.
\end{thm}

The dichotomy in the sequential behavior
of $W^{(N)}$ and $U^{(N)}$ is captured in the following proposition.
\begin{prop} \label{prop:  Strong LLN for U}
${ U^{(N)} }/{ C \log N}$
converges a.s. to $1$, while
$$\limsup_{N\to\infty} \frac{ W^{(N)} }{ C \log
N}\geq \frac{3}{2}\,.$$
In particular, $W^{(N)}/C\log N$
\emph{does not} converge a.s. to $1$.
\end{prop}

The structure of the note is as follows. In the next section, we introduce
dependency graphs and the Arratia-Goldstein-Gordon version
of the Chen-Stein method,
and perform preliminary computations. After these preliminary
computations are in place, the short Section \ref{sec-3} is devoted
to the proof of Theorem
\ref{thm: limit thm}.
Section
\ref{sec-4} is devoted to the proof of
Proposition
\ref{prop:  Strong LLN for U}.

\section{Preliminaries and auxilliary computations}
We introduce the notion of dependency graphs, and the method of
Chen and Stein to prove Poisson convergence, that will play an
important role in our proof.

\subsection{Dependency Graphs}

Let $X_1,X_2,\ldots, X_N$ be $N$ random variables.  Let $G$ be a
graph with vertices $1,2,\ldots,N$.  We use the notation
$i \sim j$ to denote two vertices connected by an edge. As
$X_i$ is not independent of itself, we define $i \sim i$ for
all $i$ (this can be thought of as requiring $G$ to have a self
loop at each vertex). $G$ is called a \emph{dependency graph} of
$\set{X_i}_{i=1}^N$ if for any vertex $i$,
\begin{equation} \label{eqn:  dependency graph}
X_i \textrm{ is independent of the set } \set{ X_j \ : \ j
\not\sim i } .
\end{equation}

The notion of dependency graphs has been introduced in connection
with the Lov\'{a}sc Local Lemma, see \cite{ProbMethod}, Chapter 5.
Some other results concerning dependency graphs are \cite{GY},
\cite{Janson}. We emphasize that  there can be many
dependency graphs associated to a collection of random variables
$\{X_i\}_{i=1}^N$.

We define two quantities associated with a dependency graph $G$
of $\{X_i\}_{i=1}^N$.
\begin{equation} \label{eqn:  B1}
B_1 = B_1(G) = \sum_{i=1}^N \sum_{j : X_j \sim X_i} \E \br{ X_i }
\E \br{ X_j } ,
\end{equation}
\begin{equation} \label{eqn:  B2}
B_2 = B_2(G) = \sum_{i=1}^N \sum_{j \neq i : X_j \sim X_i } \E
\br{ X_i X_j } .
\end{equation}

The following is a simplified version of Theorem 1 in \cite{AGG},
which in turn is an effective way to apply the Chen-Stein method:

\begin{thm}[Arratia, Goldstein, Gordon] \label{thm: AGG Thm}
Let $\set{X_i}_{i=1}^{N}$ be $N$ Bernoulli random variables with
$p_i = \E \br{ X_i } > 0$.  Set
$$ S_N = \sum_{i=1}^N X_i, \qquad \textrm{ and } \qquad \
\lambda = \E \br{ S_N} = \sum_{i=1}^N p_i . $$ %
Let $G$ be a dependency graph of $\set{X_i}_{i=1}^N$, and define
$B_1$ and $B_2$ as in (\ref{eqn:  B1}) and (\ref{eqn:  B2}).

Let $Z$ be a Poisson random variable with mean $\E \br{ Z } =
\lambda$.  Then, for any $A \subset \N$,
$$ \abs{ \Pr \br{ S_N \in A } - \Pr \br{ Z \in A} } \leq B_1 + B_2 .
$$
\end{thm}

Theorem \ref{thm: AGG Thm} is useful in proving convergence of
sums of ``almost'' independent variables to the Poisson distribution.

%\Rem Theorem \ref{thm: AGG Thm} is still true if we replace
%(\ref{eqn:  dependency graph}) with the (weaker) condition that
%for any vertex $X_i$,
%$$ \E \Br{ X_i }{ X_j \ : \ j \not\sim i } = \E \br{ X_i } . $$

%The subsequent sections are devoted to proving Proposition
%\ref{prop:  LLN} and Theorem \ref{thm: limit thm}.

\subsection{Auxilliary Calculations}

Recall that $C = {2}/{\log 2}$. Fix $\eps>0$  and set
$M = \lfloor (C+\eps) \log N
\rfloor$.
Define
\begin{equation} \label{eqn:  W', eps}
W'_{s,p} = W'^{(N)}_{s,p} \eqdef \max \Set{ 1 \leq k \leq M
 }{ \x_s = 0 \ , \ \prod_{i=1}^k \x_{s+i p {\pmod N} } = 1 }
,
\end{equation}
and $W' = \max_{s, 1 \leq p \leq N/2} W'_{s,p}$.  That is, we take truncated
versions of $W_{s,p}$ and $W$.

For $1 \leq s,p \leq N$ and $x \in \R$ define
$$ I_{s,p}(x) \eqdef \1{ W'_{s,p} > C \log N + x } , $$
and set
$$ S(x) \eqdef \sum_{\substack{ 1 \leq s \leq N \\ 
1 \leq p \leq N/2 } }  I_{s,p}(x) . $$
Note that $W' > C \log N + x$ iff $S(x) > 0$.
%To simplify the notation set $M = \lfloor (C+\eps) \log N
%\rfloor$, for $\eps$ as in (\ref{eqn:  W', eps}).
For $s,p$, let
$A(s,p) = \set{s+ip}_{i=0}^M$ be the arithmetic progression
corresponding to $I_{s,p}$.
(Recall the remark following the definition of $W^{(N)}$
that explains the restriction $p \leq N/2$.)

Let $G$ be the graph with vertex set $\set{ (s,p) }_{s \leq N ,p\leq N/2}$,
and edges defined by the relations
$$ (s,p) \sim (t,q) \quad \Longleftrightarrow \quad A(s,p)
\cap A(t,q) \neq \emptyset . $$

Fix $x \in \R$ such that $x < \eps \log N$ (for large enough $N$
this is always possible). Note that $I_{s,p}(x)$ is independent of
$\set{\x_{j {\pmod N}} \ : \ j \not\in A(s,p)}$.  Thus, $G$ is a
dependency graph of $\set{I_{s,p}(x)}_{s \leq N , p \leq N/2}$.

Define $\D_{s,p}(k)$ to be the number of pairs $t,q$ with $q \neq
p$ such that $\abs{ A(s,p) \cap A(t,q) } = k$.

The following combinatorial proposition proves to be useful.

\begin{prop} \label{prop:  combinatorial proposition for D(s,p,k)}
For all $s,p$ the following holds:
$$ \D_{s,p}(k) \leq \left\{ \begin{array}{lr}
        (M+1)^2 N & k = 1 \\
        (M+1)^2 M^2 & 2 \leq k \leq \frac{M}{2}+1 \\
        0 & k > \frac{M}{2} + 1
            \end{array} \right. $$
\end{prop}

\begin{proof}
Fix $1 \leq s,p \leq N$.

%Let $1 \leq t,q \leq N$, such that $q \neq p$.  If $\abs{ A(s,p)
%\cap A(t,q) } \geq M+1$, then $t=s$ and $q=p$, so $\D_{s,p}(k) =
%0$ for $k > M$.

Let $k \geq 2$.  Assume that $A(s,p) \cap A(t,q) = \set{x_1 < x_2
< \cdots < x_k}$.  Let $L = \mathrm{lcm}(p,q) \eqdef \min \Set{ L
}{ \exists \ a,b \in \Z \ : \ L = a p = b q }$.

\Claim $\set{x_i}_{i=1}^k$ is an arithmetic progression with
$x_{i+1} - x_i = L$.

\begin{proof}
Assume that $L = ap = bq$, for $a \neq b \geq 1$. Fix $1 \leq i
\leq k-1$. Since $x_{i+1} > x_i$ are both in $A(s,p) \cap A(t,q)$,
we get that $x_{i+1} - x_i = a' p = b' q$ for nonnegative integers
$a' \neq b' \geq 1$.  So $x_{i+1} - x_i \geq L$.

Consider $x_i + L$.  Since $x_i \in A(s,p) \cap A(t,q)$, and $L =
ap = bq$, and since $x_i < x_i + L \leq x_{i+1}$, it follows that
$x_i + L \in A(s,p) \cap A(t,q)$. So $x_i +L \geq x_{i+1}$,
concluding the proof of the claim.
\end{proof}

Let $a \neq b \geq 1$ be such that $L = ap = bq =
\mathrm{lcm}(p,q)$.  We have the following constraints:
$$ s+ (k-1) a p \leq x_1 + (k-1) L \leq s + M p \quad \textrm{ and }
\quad t+ (k-1) b q \leq x_1 + (k-1) L \leq t + M q . $$ %
Thus, $2 \leq \max \set{a,b} \leq \frac{M}{k-1}$, or: $k \leq
\frac{M}{2} + 1$.

So for $k > \frac{M}{2} + 1$ we get that $\D_{s,p}(k) = 0$.

Consider $k \leq \frac{M}{2}+1$. Since there are at most
$\frac{M}{k-1}$ choices for $a$ and for $b$, and since a choice of
$a,b$ determines $q$, we have at most $\frac{M^2}{(k-1)^2}$
choices for $q$.

\Rem This can be improved to $\frac{2}{k+1} \cdot
\frac{M^2}{(k-1)^2}$, with a slightly more careful analysis.  We
will not need this improvement.

Since $t = x_1 - i q = s + j p - iq$ for some $0 \leq i,j \leq M$,
there are at most $(M+1)^2$ choices for $t$, once we have fixed
$q$.

Thus, altogether, for $2 \leq k \leq \frac{M}{2} +1$,
$$ \D_{s,p}(k) \leq \frac{(M+1)^2 M^2}{(k-1)^2} \leq (M+1)^2 M^2 . $$

If $\abs{ A(s,p) \cap A(t,q) } = 1$ then there are at most $N$
choices for $q$ and $(M+1)^2$ choices for $t$, so $\D_{s,p}(1)
\leq (M+1)^2 N$.
\end{proof}

%We now use Proposition \ref{prop:  combinatorial proposition for
%D(s,p,k)} to prove a limit theorem for $W^{(N)}$.

Recall $G$ defined above, a dependency graph of $\set{ I_{s,p}(x)
}_{s \leq N, p \leq N/2}$.  Set
$$ B_1 = B_1(x,G) = \sum_{s \leq N , p \leq N/2} \sum_{\substack{ t,q \\ I_{t,q} \sim I_{s,p} } } \E
\br{ I_{s,p} } \E \br{ I_{t,q} } , $$ %
as in (\ref{eqn:  B1}). Also, set
$$ B_2 = B_2(x,G) = \sum_{s \leq N , p \leq N/2} \sum_{\substack{ (s,p) \neq (t,q) \\ I_{t,q}
\sim I_{s,p} } } \E \br{ I_{s,p} I_{t,q} } , $$ %
as in (\ref{eqn: B2}).

\begin{prop} \label{prop: B1(x) and B2(x) are small}
For any $\delta > 0$,
%and all $x \in \R$ such that $x < \eps \log
%N$, we have
$$ \sup_{x\in (-\infty, \eps \log N)}
B_1(x,G) + B_2(x,G) = O ( N^{\delta-1} ) . $$
%where the estimate is uniform in $x$.
\end{prop}

\begin{proof}
We have that $\E \br{ I_{t,q} } \leq 2^{- (C \log N + x + 1) }$,
for all $t,q$.

Fix $s,p$.
There is at most one value of $t$ such that $\abs{A(s,p) \cap
A(t,p)} = k$. Hence, the number of pairs $t,q$ such that
$\abs{A(s,p) \cap A(t,q)} = k$ is at most $\D_{s,p}(k) + 1$

Thus,
\begin{eqnarray*}
    B_1 & \leq & \sum_{s,p} \sum_{k=1}^{M+1} (\D_{s,p}(k) + 1) 2^{- 2 (C \log N
    + x + 1) } \\
    & \leq & 2^{-2(x+1)} \cdot \frac{1}{N^4} \sum_{s,p} \sr{ (M+1)^2 N +
    1 + \sum_{2 \leq k \leq \frac{M}{2}+1 } ((M+1)^2 M^2 + 1) } \\
    & = & O \sr{ \frac{M^2 N + M^5}{N^2} } = O \sr{ N^{\delta-1} } ,
\end{eqnarray*}
for all $\delta > 0$.

For $s,p$ and $t,q$ such that $\abs{ A(s,p) \cap A(t,q) } = k$ we
have $\E \br{ I_{s,p} I_{t,q} } \leq 2^{- 2(C \log N + x + 1) +
k}$.  Also, if $q = p$ and $A(s,p) \cap A(t,q) \neq \emptyset$,
then either $t \in A(s,p)$ or $s \in A(t,q)$.  Thus, if $t \neq
s$,
$$ \E \br{ I_{s,p} I_{t,p} } \leq \Pr \br{ \xi_s \xi_t = 0 \ , \ \xi_s
\xi_t = 1
} = 0 . $$ %
Hence,
\begin{eqnarray*}
     B_2 & \leq & \sum_{s,p} \sum_{k=1}^{M} \D_{s,p}(k) 2^{- 2(C \log N
    + x + 1) + k} \\
    & \leq & 2^{-2(x+1)} \cdot \frac{1}{N^4} \sum_{s,p} \sr{
    2 (M+1)^2 N + (M+1)^2 M^2 \cdot
    \sum_{2 \leq k \leq \frac{M}{2} +1 } 2^k } \\
    & = & O \sr{ \frac{M^2 N + M^4 2^{M/2}}{N^2} } = O \sr{
    N^{\delta - 1} } ,
\end{eqnarray*}
for all $\delta > 0$.
\end{proof}

\section{Arithmetic Progressions: Proof of Theorem
\ref{thm: limit thm}
}
\label{sec-3}

Since the proofs are very similar, we only consider the slightly harder
$W^{(N)}$.
    We write $W$ for $W^{(N)}$
whenever no confusion can occur.

%\begin{proof}[Proposition \ref{prop:  LLN}]
We begin with the following lemma:
\begin{lem}
    \label{lem-tog}
    The sequence
    ${W^{(N)}}/{C \log N}$ converges to $1$ in probability; i.e.
    for any $\delta >0$,
    $$ \lim_{N \to \infty}
    \Pr \br{ \abs{ \frac{W^{(N)}}{C \log N} - 1 } > \delta
    } = 0 . $$
Further, the convergence is almost sure on the subsequence $N_k=2^k$.
Finally, the statements hold with $U^{(N)}$ replacing $W^{(N)}$.
 %   \item Moreover, for $Y_k = W^{(2^k)} / C \log (2^k)$, we have
 %   that $Y_k$ converges a.s. to $1$.
%\end{itemize}
\end{lem}
\begin{proof}[Proof of Lemma \ref{lem-tog}]
    Again, we consider only $W^{(N)}$.
Fix
$\eps > 0$.
Note that
$$ \Pr \br{ W_{s, p} > (C+\eps) \log N } \leq 2^{-(C+\eps) \log N - 1} . $$
Thus,
\begin{equation} \label{eqn:  Prob W greater C+eps log N}
\Pr \br{ W > (C+\eps) \log N } \leq N^2 \cdot 2^{- (C+\eps) \log N
- 1} = \frac{1}{2 N^{\eps}} \longrightarrow 0 .
\end{equation}

Now let $x = - \eps \log N$, and let $Z(x)$ be a Poisson random
variable with mean
$$ \E \br{ Z(x) } = \lambda(x) = \E \br{ S(x) } = N \cdot \lfloor N/2 \rfloor  \cdot 2^{-
\lfloor (C \log N + x + 2) \rfloor} \geq 2^{\eps \log N - 3} \cdot (1-N^{-1}) . $$ %

Note that $\set{ W \leq (C-\eps) \log N }$ implies that $\set{ W'
\leq (C-\eps) \log N}$, so using Theorem \ref{thm: AGG Thm}
%OO
and Proposition \ref{prop: B1(x) and B2(x) are small},
\begin{eqnarray}
    & & \Pr \br{ W \leq (C-\eps) \log N } \leq \Pr \br{ S(x) = 0 }
    \nonumber \\
    & \leq & B_1(x,G) + B_2(x,G) + \Pr \br{ Z(x) = 0 } \nonumber \\
    & \leq & 2^{-2(x+1)} \cdot \frac{\log^5 N}{N} + \exp \sr{ - 2^{\eps
    \log N - 3} \cdot(1-N^{-1})  }  \longrightarrow 0 , \label{eqn:  Prob W less  C-eps log N}
\end{eqnarray}
for $\eps < \frac{1}{2 \log 2}$.

So for any positive $\delta < \frac{1}{4}$, we get from
(\ref{eqn:  Prob W greater C+eps log N}) and
    (\ref{eqn:  Prob W less  C-eps log N}) that
$$ \lim_{N \to \infty}
\Pr \br{ \abs{ \frac{W^{(N)}}{C \log N} - 1 } > \delta
} = 0 . $$
Further,  from the same estimates one has that with
%A quick inspection of (\ref{eqn:  Prob W greater C+eps log N}) and
%(\ref{eqn:  Prob W less  C-eps log N}) shows that for
$Y_k =
W^{(2^k)} / C \log (2^k)$, for any positive $\delta < \frac{1}{4}$,
$$ \sum_{k=1}^{\infty} \Pr \br{ \abs{ Y_k - 1} > \delta } < \infty
. $$ %
One then deduces from
the Borel-Cantelli lemma the claimed almost sure convergence.
%, we can deduce that
%$$ \Pr \br{ \limsup_{k \to \infty} \abs{ Y_k - 1 } > 0 } \leq
%\sum_{\ell=1}^{\infty} \Pr \br{ \limsup_{k \to \infty} \abs{ Y_k -
%1 } > \frac{1}{5 \ell} } = 0 , $$ proving the second assertion.
\end{proof}

% ------------------

\begin{proof}[Proof of Theorem \ref{thm: limit thm}]
As in the proof of Lemma
\ref{lem-tog},
for $x \in \R$, let $Z(x)$ be a Poisson random variable with mean
$$ \E \br{ Z(x) } = \lambda(x) = \E \br{ S(x) } = N \cdot \lfloor N/2 \rfloor \cdot 2^{-
\lfloor (C \log N + x + 2) \rfloor} . $$ %
If $C\log N +x \in \Z$, then $\lambda(x) = 2^{-(x+2)} \cdot \frac{\lfloor N/2 \rfloor }{N}$.

Note that $W' > C \log N + x$ iff $S(x) > 0$. By Theorem \ref{thm:
AGG Thm}
%OO
and Proposition \ref{prop: B1(x) and B2(x) are small},
\begin{eqnarray*}
    & & \abs{ \Pr \br{ W' > C \log N + x } - \Pr \br{ Z(x) \neq 0 }
    } = \abs{ \Pr \br{ S(x) > 0 } - \Pr \br{ Z(x) > 0 } } \\
    & \leq & B_1(x,G) + B_2(x,G) = O \sr{ N^{\delta-1} } .
\end{eqnarray*}
We also have the equality
$$ \set{ W > C \log N + x } = \set{ W > (C+\eps) \log N } \bigcup
\set{ W' > C \log N + x } . $$ %
Thus, for $0 < \delta < 1$,
\begin{eqnarray*}
    & & \abs{ \Pr \br{ W \leq C \log N + x } - e^{-\lambda(x)} }
    \\
    & \leq & \Pr \br{ W > (C+\eps) \log N } + \abs{ \Pr \br{ W' > C \log N + x } - \Pr \br{ Z(x) \neq 0
    } } \\
    & \leq & O \sr{ N^{-\eps} } + O \sr{ N^{\delta-1} } .
\end{eqnarray*}

Let $\set{x_N}$ be a sequence such that $C \log N + x_N \in \Z$
for all $N$. If $\inf_N x_N \geq b \in \R$, then $\exp \sr{
\lambda(x_N) }$ is a bounded sequence.  Thus,
$$ \abs{ \exp \sr{ \lambda(x_N) } \Pr \br{ W^{(N)} \leq C \log N + x } -
1 } = O \sr{ N^{-\eps} } \longrightarrow 0 . $$
Since 
%$\lambda(x_N) = 2^{-(x_N+3)}$ if $N$ is even and
%if $N$ is odd 
$$ \abs{ \exp ( \lambda(x_N)) - \exp (2^{-(x+3)}) } \leq \exp (2^{-b-3}) \cdot \sr{
\exp ( 2^{-b-3} N^{-1}) - 1 } = O(N^{-1}) , $$
we have the theorem (we have used again that $\inf_N x_N \geq b \in \R$).
\end{proof}

\section{Convergence in Probability vs. a.s. Convergence}
\label{sec-4}
%
%In Proposition \ref{prop:  LLN} we proved that the sequence
%$\frac{ W^{(N)} }{\log N}$ converges to $\frac{2}{\log 2}$ in
%probability.  The following question immediately arises:  Does
%a.s. convergence also hold? This question is strengthened by the
%following observations:
%
%In considering $W^{(N)}$ we considered arithmetic progressions
%$\pmod N$.  We can also consider arithmetic progressions in
%$\set{1,\ldots,N}$ that do not wrap around the torus; i.e., define
%$$ U_{s,p} = U_{s,p}^{(N)} \eqdef \max \Set{ 1 \leq k \leq N }{ \x_s = 0 \ , \
%\prod_{i=1}^k \x_{s+i p} = 1 } , $$ %
%and
%$$ U = U^{(N)} \eqdef \max_{s, p} U_{s,p} , $$
%where for we only consider $s,p,k$ such that
%$\Set{s+ip}{i=0,1,\ldots,k} \subseteq \set{1,\ldots,N}$.
%
%The following proposition can be proved similarly to the way
%Proposition \ref{prop:  LLN} is proved.
%
%\begin{prop}
%For $C = \frac{2}{\log 2}$, $\frac{U^{(N)}}{C \log N}$ converges
%in probability to $1$.  Moreover, for $Y_k = U^{(2^k)} / C \log
%(2^k)$, we have that $Y_k$ converges a.s. to $1$.
%\end{prop}
%
%
We begin with the following easy consequence of Lemma
\ref{lem-tog} applied to $U^{(N)}$.
%This proposition is enough to prove a.s. convergence of the whole
%sequence.
%
\begin{prop} \label{prop:  SStrong LLN for U}
%For $C = \frac{2}{\log 2}$,
${ U^{(N)} }/{ C \log N}$
converges a.s. to $1$.
\end{prop}

\begin{proof}
The main observation is that $U^{(N)}$ is a monotone increasing sequence.
That is, a.s. for all $N$, $U^{(N)} \leq U^{(N+1)}$. Thus, a.s.
for all $N$, setting $k = \lfloor \log_2 N \rfloor$, we have
\begin{eqnarray*}
    \frac{ U^{(2^k)}}{ \log (2^k)} \cdot \frac{ \log (2^k) }{ \log (2^{k+1}) }
    & \leq & \frac{U^{(N)} }{ \log N } \leq \frac{ U^{(2^{k+1})} }{ \log (2^{k+1}) } \cdot \frac{
    \log (2^{k+1}) }{ \log (2^k) }
\end{eqnarray*}

Since
%OO
${ U^{(2^k)}}/{ \log (2^k)}$ converges a.s. to $C$ and
${ \log (2^k) }/{ \log (2^{k+1}) }$ converges a.s. to $1$, we
get a.s. convergence of $\frac{U^{(N)} }{ \log N }$ to $C$.
\end{proof}

We turn to the
\begin{proof}[Proof of Proposition  \ref{prop:  Strong LLN for U}]
In view of Proposition  \ref{prop:  SStrong LLN for U}, it remains only
to consider the statement concerning $W^{(N)}$. Toward this end,
%
%
%Recall that $W^{(N)}$ is the same as $U^{(N)}$, only taking into
%consideration also arithmetic progressions $\pmod N$; i.e. instead
%%of considering the interval from $1$ to $N$, considering the
%torus.  This seemingly insignificant perturbation, surprisingly
%changes the behavior of the longest arithmetic progression.
%
%
%
%\begin{prop} \label{prop:  No a.s. convergence}
%For $C = \frac{2}{\log 2}$, the sequence $\frac{ W^{(N)} }{ C \log
%N}$ \emph{does not} converge a.s. to $1$.
%\end{prop}
%
%
%\begin{proof}
fix $0 < \beta < 1$.  Let $M_N = \lceil (2+\beta) \log_2 N
\rceil$. So $N^{2+\beta} \leq 2^{M_N} \leq 2 N^{2+\beta}$.  Define
$$ I(s,p,N) = \1{ W_{s,p}^{(N)} \geq M_N } . $$
That is, $I(s,p,N)$ is the indicator function of the event that
$\x_s = 0$ and $\prod_{i=1}^{M_N} \x_{s+ip {\pmod N} }=1$.
Set
$${\rm Cov}( (s,p,N), (s',p',N'))=
    \E \br{ I(s,p,N) I(s',p',N') } - \E \br{ I(s,p,N) } \E \br{
    I(s',p',N') } .$$

For simplicity of notation, for $a,b \in \Z$, we denote $[a,b] =
\Set{ z \in \Z }{ a \leq z \leq b }$.
Let ${\cal L}_{n}=\{(s,p,N): N\in[n,2n], s\in[1,N], p\in[N/6,N/5]\}$.
\begin{lem} \label{lem:  bound on second moment}
The following holds:
%For large enough $n$ the following holds:  Let $n \leq N \leq 2n$
%and let $1 \leq s,p \leq N$. Then,
%\begin{eqnarray}
$$     \sum_{(s,p,N)\in {\cal L}_n}
     \sum_{(s',p',N')\in {\cal L}_n}
%N,N' = n}^{2n} \sum_{s,p =1}^N \sum_{s',p' =1}^{N'}
{\rm Cov}( (s,p,N), (s',p',N'))\leq
    %\E \br{ I(s,p,N) I(s',p',N') } - \E \br{ I(s,p,N) } \E \br{
    %I(s',p',N') } \nonumber \\
    %& \leq &
    O \sr { n^{1-\beta}\log^7(n)   }.
    $$
%\end{eqnarray}
\end{lem}

\begin{proof}
For any
$N,N' \in [n,2n]$,
$s,p \in [1,N]$ and $s',p' \in [1,N']$, we have
that
\begin{equation}
    \label{eq-180707a}
    {\rm Cov}((s,p,N), (s',p',N') )
 \leq 2^{-(M_N + M_{N'})} \sr{ 2^k - 1},
 \end{equation}
with
$$ k = \abs{ \set{ s+ip {\pmod N} }_{ i \in [1,M_N] } \cap \set{ s' +
j p' {\pmod {N'}} }_{ j \in [1,M_{N'}] } } . $$ %
Thus, the proof of
the lemma is based on
    controlling the cardinality of the collection of triples
    $(s',p',N')\in {\cal L}_n$, whose associated
    arithmetic progression intersects in a prescribed number
    of points the arithmetic progression associated with
    a given triple $(s,p,N)\in {\cal L}_n$. We divide our estimates
    into three: intersection at one point, intersection at two
    points or more, and intersection at $2C\log (2n)/5$ points or more.

For $N,N' \in [n,2n]$ and $p \in [1,N]$ define $\T(N,N',p)$ to be
the set of all triples $(s,s',p')$ such that $s \in [1,N]$, $s',p'
\in [1,N']$ and
\begin{equation} \label{eqn:  dfn of T}
\abs{ \set{ s+ip {\pmod N} }_{ i \in [1,M_N] } \cap \set{ s'+jp'
{\pmod {N'}} }_{ j \in [1,M_{N'}] } } = 1 .
\end{equation}
Similarly, define $\S(N,N',p)$ to be
the set of all such triples $(s,s',p')$ such that
%$s \in [1,N]$, $s',p'
%\in [1,N']$ and
\begin{equation} \label{eqn:  dfn of S}
\abs{ \set{ s+ip {\pmod N} }_{ i \in [1,M_N] } \cap \set{ s'+jp'
{\pmod {N'}} }_{ j \in [1,M_{N'}] } } \geq 2 .
\end{equation}
Finally, define $\U(s,p,N)$ the set of all triples
$(s',p',N')\in {\cal L}_n$ such that
\begin{equation} \label{eqn:  dfn of U}
\abs{ \set{ s+ip {\pmod N} }_{ i \in [1,M_N] } \cap \set{ s'+jp'
{\pmod {N'}} }_{ j \in [1,M_{N'}] } } \geq  2M_{N'}/5 .
\end{equation}
We have the following estimates.

\begin{prop} \label{prop:  bound on T}
For large enough $n$, the following holds: for all $N,N' \in
[n,2n]$ and $p \in [1,N]$,
$$ \abs{ \T(N,N',p) } \leq n^2 \log^5(n) . $$
\end{prop}

\begin{prop} \label{prop:  bound on S}
For large enough $n$, the following holds: for all $N,N' \in
[n,2n]$ and $p \in [1,N]$,
$$ \abs{ \S(N,N',p) } \leq n \log^9(n) . $$
\end{prop}

\begin{prop} \label{prop:  bound on U}
    For large enough $n$,
    the following holds: for all $(s,p,N)\in {\cal L}_n$,
$$ \abs{ \U(s,p,N) } \leq  \log^7(n) . $$
\end{prop}
Assuming Propositions
\ref{prop:  bound on T},
\ref{prop:  bound on S},
\ref{prop:  bound on U}, we have
%conclude the proof of the lemma.
%Note that
\begin{eqnarray} \label{eqn:  sum of s,s',p'}
    & &
    \sum_{(s,p,N)\in {\cal L}_n}
    \sum_{(s',p',N')\in {\cal L}_n}
 {\rm Cov}((s,p,N), (s',p',N') ) \\
 &\leq &   \sum_{N,N'=n}^{2N}
 \sum_{p=1}^N
 2^{-(M_N+M_{N'})}
 \abs{\T(N,N',p)}
 + \sum_{N,N'=n}^{2n}\sum_{p=1}^N
 2^{-(M_N+M_{N'})}
 %OO
 \abs{\S(N,N',p)} 2^{2M_{N'}/5}  \nonumber \\
 &&+
    \sum_{(s,p,N)\in {\cal L}_n}
 \abs{\U(s,p,N)} 2^{-M_n}  \nonumber \\
 & \leq & O(n^{1-2\beta} \log^5(n))+
 O(n^{4/5-8\beta/5}\log^9(n))+O(n^{1-\beta}\log^7(n))
 \leq O(n^{1-\beta}\log^7(n))\,,
 \nonumber
 %2^{-(M_N + M_{N'})} n^2 \log^5(n) + n2^{-M_{N'}} \log^9(n) .
\end{eqnarray}
which completes the proof of the lemma.
\end{proof}

Returning to the proof of Proposition
\ref{prop:  Strong LLN for U},
let
$$ \Lambda(n) = \sum_{(s,p,N)\in {\cal L}_n}
 I(s,p,N)  $$
and
note that for all large enough $n$,
\begin{eqnarray} \label{eqn:  first moment lower bound}
    \E[\Lambda(n)]=\sum_{(s,p,N)\in {\cal L}_n}
\E \br{ I(s,p,N) }
    & = & \Omega \sr{  n^{1-\beta}  }
\end{eqnarray}
while from Lemma \ref{lem:  bound on second moment},
\begin{eqnarray} \label{eqn:var}
    \mbox{\rm Var}\sr{\Lambda(n)}=
    \mbox{\rm Var}\sr{\sum_{(s,p,N)\in {\cal L}_n}
    \br{ I(s,p,N) }}
    & = & O \sr{  n^{1-\beta}\log^7(n)  }\,.
\end{eqnarray}
Thus,
\begin{equation}
    \label{eq-180707b}
    \Pr\br{\Lambda(n)=0}\leq \frac{\mbox{\rm Var}\sr{\Lambda(n)}}{
    \sr{\E \Lambda(n)}^2}=O\sr{n^{\beta-1}\log^7(n)}\,.
\end{equation}
Since
$$ \Pr \br{ \exists \ N \in [n,2n] \ : \ W^{(N)} >
\frac{2+\beta}{\log 2} \log N } \geq \Pr \br{ \Lambda(n) > 0 } ,
$$
it follows from (\ref{eq-180707b}) that
%
%Now, by Cauchy-Schwartz,
%$$ \E \br{ \Lambda(n) } = \E \br{ \Lambda(n) \1{ \Lambda(n) > 0} }
%\leq \sqrt{ \E \br{ (\Lambda(n) )^2 } \Pr \br{ \Lambda(n) > 0 } }
%. $$ %
%By (\ref{eqn: first moment lower bound}),
%$$ \E \br{ \Lambda(n) } = \Omega \sr{ n^{1-\beta} } . $$
%By Lemma \ref{lem:  bound on second moment} and (\ref{eqn:  first
%moment lower bound}),
%$$ \E \br{ \Lambda(n)^2 } = O \sr{ n^{2-2\beta} \cdot \br{ 1 +
%\log^9(n) n^{\beta-1} } } . $$ %
%Thus,
%$$ \Pr \br{ \exists \ N \in [n,2n] \ : \ W^{(N)} >
%\frac{2+\beta}{\log 2} \log N } \geq  1 - O (\log^9(n) n^{\beta-1}) . $$ %
%
%Taking $n=2^k$, we get that
$$
\sum_{k=0}^{\infty} \Pr
\br{ \max_{N \in [2^k, 2^{k+1}]} \frac{W^{(N)} }{\log N}
\leq \frac{2+\beta}{\log 2} } < \infty . $$
By the
%OO
Borel-Cantelli lemma, we get that a.s.
$$ \limsup_{N \to \infty} \frac{W^{(N)} }{ \log N} \geq
\limsup_{k \to \infty} \max_{N \in [2^k, 2^{k+1}]} \frac{W^{(N)}
}{\log N} > \frac{2+\beta}{\log 2} . $$ %
Since $\beta\in (0,1)$ is arbitrary, this completes
the proof of Proposition
\ref{prop:  Strong LLN for U}.
\end{proof}

\begin{proof}[Proof of Proposition  \ref{prop:  bound on T}]
If $(s,s',p') \in \T(N,N',p)$ then (\ref{eqn:  dfn of T}) implies
that there exist $i \in [1,M_N]$, $j \in [1, M_{N'}]$, $k_i \in
[0,M_N]$ and $k'_j \in [0,M_{N'}]$ such that
\begin{equation}
    s + ip - k_i N = s' + j p' - k'_j N'
\end{equation}
There are at most $(2n)^2$ choices for $s'$ and $p'$.  There
exists some universal constant $K$ such that there are at most $K
\log(n)$ choices for each of $i,j,k_i,k'_j$. Choosing
$s',p',i,j,k_i,k'_j$ determines $s$. Thus, we have shown that
$\abs{ \T(N,N',p) } \leq 4 K n^2 \log^4(n) \leq n^2 \log^5 (n)$
for large enough $n$.
\end{proof}
\begin{proof}[Proof of Proposition  \ref{prop:  bound on S}]
If $(s,s',p') \in \S(N,N',p)$ then (\ref{eqn:  dfn of S}) implies
that there exist $i,r \in[1,M_N]$ and $j,\ell \in [1,M_{N'}]$ such
that $(i,j)\neq (r,\ell)$ and
\begin{eqnarray} \label{eqn:  s+ip mod N = s'+j p' mod N'}
    & s + ip {\pmod N} = s'+jp' {\pmod {N'}} &  \\
    \textrm{ and } & s + r p {\pmod N} = s'+\ell p' {\pmod {N'}} .
    & \nonumber
\end{eqnarray}
Note that for any $i \in [1,M_N]$ there exists $k_i \in [0,M_N]$
such that $s+ip {\pmod N} = s+ip - k_i N$.  Similarly, for any $j
\in [1,M_{N'}]$ there exists $k'_j \in [0,M_{N'}]$ such that $s'+j
p' {\pmod {N'}} = s' + j p' - k'_j N'$.  Plugging this into
(\ref{eqn: s+ip mod N = s'+j p' mod N'}), and subtracting
equations, we get that there exist $i,r \in [1,M_N]$, $j,\ell \in
[1,M_{N'}]$, $k_i,k_r \in [0,M_N]$ and $k'_j, k'_{\ell} \in
[0,M_{N'}]$ such that
\begin{equation} \label{eqn:  rp + aN = (j-ell)p' + b N'}
(r-i) p + (k_i - k_{r}) N = (\ell-j) p' + (k'_j - k'_{\ell}) N'
.
\end{equation}
There exists some universal constant $K$ such that there are at
most $K \log (n)$ choices for each of
$i,r,j,\ell,k_i,k_r,k'_j,k'_\ell$, and $2n$ choices for $s$. After choosing
$i,r,j,\ell,k_i,k_r,k'_j,k'_\ell,s$, (\ref{eqn: s+ip mod N = s'+j p'
mod N'}) and (\ref{eqn: rp + aN = (j-ell)p' + b N'}) determine
$s'$ and $p'$.  Thus, we have shown that
for large enough $n$,
$$\abs{\S(N,N',p) } \leq
2n\sr{ K \log(n) }^8 \leq n\log^9 (n).$$
\end{proof}

\begin{proof}[Proof of Proposition  \ref{prop:  bound on U}]
Let $A={ \Set{ s+ip {\pmod N} }{ i \in [1,M_N] } }$, and let
$(s',p',N')\in \U(s,p,N)$.
 For $i = 1,6,11,\ldots$, let
 $$ Z_i = \Set{ s' + (i+r) p' {\pmod {N'}} }{ r = 0,1,\ldots,4 } . $$
This is a partition of the arithmetic progression into packets of
five elements.  We then have, by the definition of
$\U(s,p,N)$,
$$ \frac{2M_{N'}}{5}  \leq \sum_i \abs{Z_i} \leq \frac{M_{N'}}{5} \max_i
\abs{Z_i} . $$ %
So there exists some set $Z_i$ such that $\abs{Z_i} \geq 2$. This
implies that there exist $x < y \in A \cap [1,N']$, $i \in [1,M_{N'}]$,
and $r \in [1,4]$ such that
\begin{eqnarray} \label{eqn:  x,y, first two points}
    & & s + i p' {\pmod {N'} } = x  \ ,  \\
    & & s + (i+r) p' {\pmod {N'}} = y  \nonumber .
\end{eqnarray}
Subtracting equations, and using the fact that $r p' < N'$, we get
that
\begin{equation} \label{eqn:  fixing p}
    r p' = y-x .
\end{equation}

Moreover, (\ref{eqn:  dfn of U}) also implies that there must
exist an integer $j$ (perhaps negative) with $\frac{1}{5} M_{N'} \leq
\abs{j} \leq M_{N'}$, and $z \in A \cap [1,N']$, such that
\begin{equation} \label{eqn:  z, third point}
    s + (i+j) p' {\pmod {N'}} = z .
\end{equation}
For large enough $n$, we have that $\abs{j} \geq 7$.  Since $7 p' >
N'$, we get by subtracting (\ref{eqn:  x,y, first two points}) from
(\ref{eqn:  z, third point}),
\begin{equation} \label{eqn:  fixing N}
    j p' + k N' = z-x ,
\end{equation}
for some $k \neq 0$, such that $\abs{k} \leq M_{N'}$.

Since $k r \neq 0$, equations (\ref{eqn:  fixing p}) and
(\ref{eqn:  fixing N}) have at most one solution for $p',N'$, in
terms of $x,y,z,r,j$ and $k$.  Since there are at most $|A|^3\leq
M_N^3$
choices for $x,y$ and $z$, at most $4$ choices for $r$, and at
most $4 M_{N'}^2$ choices for $j$ and $k$, we get that there are at
most $16 |A|^3 M_{N'}^2$ choices for $p',N'$.  Also, there are at most $M_{N'}$
choices for $i$ in (\ref{eqn:  x,y, first two points}), and fixing
$p',N',x$ and $i$ determines $s'$.  Thus,
$$ \abs{ \U( (s,p,N)) } \leq 16 M_N^3 M_{N'}^3 . $$
\end{proof}

{\bf Open Problem.} We conjecture that in fact,
$$ \limsup_{N \to \infty} \frac{W^{(N)} }{ \log N} = \frac{3}{2}.$$

%OO
{\bf Acknowledgements.}  We thank Tim Austin, Ori Gurel-Gurevich,
Gady Kozma and Eric Shellef for useful discussions.
We also thank Min-Zhi Zhao and Hui-Zeng Zhang for pointing out a computational error
in a previous version of this note.  Their paper \cite{ZZ12} 
provides extensions of our results to$p \neq 1/2$, a resolution of our
conjecture above, and several other refinements.

% -----------------------------------------------------------------------------------
% ---------------------------- BIBLIOGRAPHY -----------------------------------------
% -----------------------------------------------------------------------------------

% -----------------------------------------------------------------------------------
% ---------------------------- END OF BIBLIOGRAPHY ----------------------------------
% -----------------------------------------------------------------------------------

\end{document}